\documentclass[a4paper,10pt]{article}
\usepackage[utf8x]{inputenc}

\usepackage{amsmath}
\usepackage{amsfonts} 
\usepackage{latexsym} 
\usepackage{amssymb}
\usepackage{amsthm}   
\usepackage{mathrsfs} 

\renewcommand{\P}{\mathbb{P}}

\newcommand{\Q}{\mathbb{Q}}

\newcommand{\Z}{\mathbb{Z}}
\newcommand{\A}{\mathbb{A}}
\newcommand{\G}{\mathbb{G}}
\newcommand{\Ci}{\mathcal{C}}

\newcommand{\M}{\mathcal{M}}
\renewcommand{\O}{\mathcal{O}}

\renewcommand{\k}{k}

\renewcommand{\bar}[1]{\overline{#1}}

\renewcommand{\phi}{\varphi}
\renewcommand{\epsilon}{\varepsilon}

\DeclareMathOperator{\Pic}{Pic}

\providecommand{\abs}[1]{\left\lvert#1\right\rvert}

\newtheorem{thm}{Theorem}
\newtheorem{lemma}[thm]{Lemma}

\theoremstyle{definition}

\theoremstyle{remark}
\newtheorem{rmk}[thm]{Remark}

\title{Quadratic integral solutions to double Pell equations}
\author{Francesco Veneziano\\
f.veneziano@sns.it}

\begin{document}

\maketitle

\begin{abstract}
 We study the quadratic integral points---that is, ($S$-)integral points defined over any extension of degree two of the base field---on a curve defined in $\P_3$ by a system of two Pell equations. Such points belong to three families explicitly described, or belong to a finite set whose cardinality may be explicitly bounded in terms of the base field, the equations defining the curve and the set $S$. We exploit the peculiar geometry of the curve to adapt the proof of a theorem of Vojta, which in this case does not apply.
\end{abstract}
\section{Introduction}
Let $\Ci$ be an irreducible curve defined over a number field $\k$. We know, depending on the genus, the general structure of $\Ci(\Q)$; on the other hand, of course $\Ci(\overline{\Q})$ is always infinite, so it is natural to ask what happens if we consider algebraic points up to some fixed degree over $\Q$.

Abramovich and Harris in \cite{AbramovichHarris}, and also Silverman and Vojta in \cite{Report3,MR1085268,Report6} were among the first to study the set of points $P\in\Ci(\bar{\Q})$ such that $[\k(P):\k]\leq d$, in particular whether it is infinite or not.

In this paper we deal with the analogous problem for {\em integral} points, sticking to the special case of $d=2$.
This case was already studied in generality by Corvaja and Zannier in \cite[Corollary 1]{CZ3} using Schmidt's Subspace Theorem.

We will further specialise the problem to curves defined by a double Pell equation, such as, for example,
\begin{equation}
	\begin{cases}\label{DefiningExample}
		y^2=2x^2+1\\
		z^2=3x^2+1.
	\end{cases}
\end{equation}
Double equations of this kind are historically relevant, being among the first curves of genus 1 ever studied.\footnote{ See for instance \cite{WeilNT} or \cite{MordellDiophEq} for a discussion of some classical cases.}

Of course, by Siegel's theorem on integral points on curves, there are only finitely many solutions in $\Z$ to \eqref{DefiningExample}, but it may be easily seen that there are infinitely many solutions in algebraic integers of degree $2$ over $\Q$: in fact it is classically known that each Pell equation $x^2-dy^2=1$ has infinitely many integral solutions when $d$ is a positive integer not a square, so one can solve the first equation to find infinitely many $(x_n,y_n)$ such that $y_n^2=2x_n^2+1$, and then set $z_n=\pm\sqrt{3x_n^2+1}$; in this way we can actually find three infinite families of solutions. We will fully describe the set of quadratic integral points on these curves and give a geometrical meaning to the families just mentioned.

\bigskip

Going back to the general context, Abramovich and Harris conjectured that the set of $P\in\Ci(\bar{\Q})$ such that $[\k(P):\k]\leq d$ is infinite (up to a finite extension of the basefield $\k$) if and only if there exists some non constant morphism $\phi:\Ci\to X$ of degree at most $d$, where $X$ is either $\P_1$ or an elliptic curve $E$ with $\abs{E(\k)}=\infty$.

Note that this condition is obviously sufficient; if $\phi$ is defined over $\k$ the preimage through $\phi$ of a rational point in $X$ is a point of degree at most $d$ in $\Ci$.

While the general conjecture was proved false by Dabarre and Fahlaoui in \cite{MR1241949}, Abramovich and Harris managed to prove it in some cases as, for example, when $d=2,3$.

The first step in their proofs was to consider the $d$-fold symmetric product of the curve, $\Ci^{(d)}$; this is a variety whose set of points may be identified with the set of all unordered $d$-tuples of points of $\Ci$. Points of degree at most $d$ on $\Ci$ naturally correspond to $\k$-rational points on $\Ci^{(d)}$, so the said authors could work on the variety $\Ci^{(d)}$ and apply results by Faltings after mapping $\Ci^{(d)}$ in an abelian variety. \footnote{For an account on this problem see \cite{MR1288998}.}

If we consider only points of degree $2$ over the base field, the existence of infinitely many quadratic rational points is surely necessary to have infinitely many quadratic {\em integral} points, but it is not sufficient, as we shall see.

In \cite{CZ3} Corvaja and Zannier prove a theorem on integral points on surfaces and apply it to $\Ci^{(2)}$ to get the following theorem:  
\begin{thm}\label{ThmCZQuadratic}
Let $\tilde\Ci$ be a projective non-singular curve defined over a number field
$\k$, and $\Ci=\tilde\Ci\setminus\{Q_1,\dotsc,Q_r\}$ be an open affine subset,
for distinct $Q_i\in\tilde\Ci(\k)$. Then
\begin{enumerate}
\item If $r\geq 5$ then $\Ci$ contains only finitely-many quadratic-integral
points over $\k$;
\item If $r=4$ there exist finitely many rational maps $\psi:\tilde\Ci\to\P_1$
of degree $2$ such that all but a finite number of the quadratic-integral
points on $\Ci$ over $\k$ are sent to $\P_1(\k)$ by at least one of these
maps.
\end{enumerate}
\end{thm}
As it happens for the structure of the $\k$-integral points, what matters in this problem is the number of points at infinity.

\section{Setting and statement}
As mentioned above, this paper will study quadratic integral points on some special curves in $\P_3$ of genus 1 and with $4$ points at infinity, defined by a double Pell-like equation; let then $a,b,c,d$ be algebraic integers such that $abcd\neq 0$ and $ad-bc\neq 0$, and let $\Ci$ be the affine
curve in $\A^3$ defined by 
\begin{equation}
	\begin{cases}\label{DefiningEquation}
		y^2=ax^2+c\\
		z^2=bx^2+d.
	\end{cases}
\end{equation}
Let $\tilde{\Ci}$ be its projective completion, defined by homogeneous equations
\begin{equation*}
	\begin{cases}
		Y^2=aX^2+cW^2\\
		Z^2=bX^2+dW^2.
	\end{cases}
\end{equation*}

Let us indicate with $P_1,P_2,P_3,P_4$ the four points at infinity of
$\tilde{\Ci}\setminus \Ci$, which, in the coordinates $(X:Y:Z:W)$ are the points
\begin{align*}
	&P_1: (1 :\sqrt{a} :\sqrt{b} :0)\\
	&P_2: (1 :\sqrt{a} :-\sqrt{b} :0)\\
	&P_3: (1 :-\sqrt{a} :\sqrt{b} :0)\\
	&P_4: (1 :-\sqrt{a} :-\sqrt{b} :0).
\end{align*}
Let $\k$ be a number field containing $\sqrt{a},\sqrt{b},c$ and $d$, and let $S\subseteq\M_\k$ a finite set of
absolute values of $\k$ containing all the archimedean ones and all the primes in $cd(bc-ad)$; let $s$ be the cardinality of $S$.

\begin{thm}\label{Main}
The set of quadratic $S$-integral points on $\Ci$ is the union of:
\begin{itemize}
 \item Three families consisting of the preimages through the three maps $\Ci\to \P_1$
     \begin{align*}
	  (x,y,z)\mapsto (x,y)\\
	  (x,y,z)\mapsto (x,z)\\
	  (x,y,z)\mapsto (y,z)
     \end{align*}
  of the $S$-integral points of $\P_1$;
 \item A finite set of cardinality at most $2^{2835s+3}$;
 \item A finite and effectively computable set whose cardinality is at most $3\cdot 2^{1121(s+H-1)+1}$, where $H$ is the class number of $\O_S$.
\end{itemize}
\end{thm}
\begin{rmk}
 The three families of quadratic points are indeed easy to spot. If for example the first of the two Pell equations defining $\Ci$ has integral solutions $(x_n,y_n)$, then $(x_n,y_n,\pm\sqrt{bx_n^2+d})$ are quadratic integral points on $\Ci$.
Similarly if $(x_n,z_n)$ are integral solutions to the second equation, then we have a second family $(x_n,\pm\sqrt{ax_n^2+c},z_n)$, and if $(y_n,z_n)$ are integral solutions to $by_n^2-az_n^2=bc-ad$, we get $\left(\pm\sqrt{\frac{y_n^2-c}{a}},y_n,z_n\right)$ for a subsequence of $(y_n,z_n)$.
\end{rmk}
\begin{rmk}
We also note that this theorem allows for a bound on the number of exceptional solutions outside the three infinite families, while Theorem \ref{ThmCZQuadratic} does not.

While the Subspace theorem is not effective, there is a semi-effective version due to Evertse which provides an explicit bound for the {\em number} of exceptional hyperplanes; it is not possible, however, to use it in the proof of Theorem \ref{ThmCZQuadratic} to get such a bound for the quadratic integral points.
 This is because, in the proof of Theorem \ref{ThmCZQuadratic}, the Subspace theorem is applied to a surface, and each exceptional hyperplane for the Subspace Theorem gives an exceptional curve on this surface which may contain points corresponding to quadratic integral points; but even if we can bound the number of such curves, they are not effectively computable and so it is not possible to bound the number of points on them.
\end{rmk}
\begin{rmk}
 We finally remark that we can not expect in general a bound for the finite sets of points which is uniform in the coefficients $a,b,c,d$; for example it is known\footnote{This result is essentially due to Mahler in \cite{Mahlerx3y3c} and improved by Silverman in \cite{MR703458}.} that there is a constant $C>0$ such that for infinitely many positive integers $A$, the number of positive integer solutions to the equation $X^3+Y^3=A$ exceeds $C\sqrt[3]{\log A}$.
\end{rmk}

\section{Sketch of the proof}
To study quadratic integral points we proceed as \cite{AbramovichHarris} and \cite{CZ3} and consider the symmetric square of the affine curve $\Ci$, which is defined as the quotient of the product $\Ci\times\Ci$ by the involution which exchanges the couple $(P,Q)$  with the couple $(Q,P)$. This quotient is an irreducible surface which can be identified with the set of unordered couples of points, and the quotient map is
\begin{align*}
 \Ci\times\Ci\xrightarrow{\phi} \Ci^{(2)}\\
 (P,Q)\mapsto \{P,Q\}.
\end{align*}
Consider a point $P$ defined on a field of degree $2$ over $\k$, and let $P'$ be its conjugate. Then the point $\{P,P'\}$ on $\Ci^{(2)}$ is fixed by every Galois automorphism of $\bar{\Q}/\k$, and hence is defined over $\k$.
If furthermore $P$ is an integral point so is the point $\{P,P'\}$. \footnote{See, for example, \cite{SerreAlgebrGrpsClassFlds} for more on the symmetric product.}

The special geometry of curves defined by double Pell equations carries to their symmetric square and allows one to study quadratic integral points directly and to give an explicit description of the maps mentioned in Theorem \ref{ThmCZQuadratic}; for the proof we will mimic the proof of a theorem by Vojta (\cite{VojtaDiophAppVDT}):
\begin{thm}[Vojta]\label{ThmVojtaDiv}
 Let $V$ a projective, nonsingular variety over a number field $\k$. Let $r$ be
the rank of the group of $\k$-rational points of $\Pic^0(V)$, $\rho$ the
rank of the Néron-Severi group of $V$ and $D$ a divisor with at least $\dim
V+r+\rho+1$ distinct irreducible components, all defined over $\k$. Then all
sets of quasi-$S$-integral points on $V\setminus\abs{D}$ are degenerate.
\end{thm}
\begin{rmk}
 This theorem has been improved later by Vojta himself, removing the assumptions on the $\Pic^0$, in \cite{VojtaSubvSemiabVar} and also by Noguchi and Winkelmann in \cite{Report4}.
\end{rmk}

While this theorem does not apply to our case, the proof adapts well because on the curve $\Ci$ the difference of any two points at infinity is torsion in the Picard group.

\bigskip

The first step of the proof, following the strategy already illustrated, will be to study instead the structure of the set of integral points on the surface $\Ci^{(2)}$ obtained by taking the symmetric square of the original curve.

We will then follow the proof of Theorem \ref{ThmVojtaDiv} and build three functions $\alpha, \beta, \gamma$ without zeroes and poles on $\Ci^{(2)}$. These functions, up to constant factors, take integral points of the symmetric square to $S$-units.

We will then find a relation among $\alpha, \beta, \gamma$, so that taking them as coordinates gives a map from $\Ci^{(2)}$ to the subvariety of $\G_m^3$ defined by this relation; it will turn out that this relation is linear.

This means that the functions $\alpha, \beta, \gamma$ take integral points of $\Ci^{(2)}$ to solutions of the $S$-unit equation, which is the object of the following theorem.
\begin{thm}[Evertse, \cite{EvertseTheNumberOfSolutions}]\label{EvertseSUnitEffective}
	Let $\k$ be a number field, let $S$ be a finite set of places of $\k$
containing all the archimedean ones, let $s=|S|$ and let $a_1,\dotsc,a_n\in
\overline{\Q}^*$. Let $A'(a_1,\dotsc,a_n;\O_S^*)$ be the number of
non degenerate solutions
$x_1,\dotsc,x_n\in\O_S^*$ to the equation
	\begin{equation}\label{CoefficientSUnit}
		a_1x_1+\dotsb +a_nx_n=1,
	\end{equation}
where a solution $(x_1,\dotsc,x_n)$ is called {\em degenerate} if there is a proper vanishing subsum in the left hand side of \eqref{CoefficientSUnit}.
Then
\[
	A'(a_1,\dotsc,a_n;\O_S^*)\leq 2^{35 n^4 s}.
\]
\end{thm}
Using this theorem we will then proceed to bound the number of non degenerate solutions and to examine degenerate solutions, which come from special subvarieties.

We will show that three of these special subvarieties are curves of genus 1, hence each gives only a finite number of quadratic integral points, while the other three have genus 0 and give three families of quadratic integral points.

\medskip

The strategy outlined here would also work for any variety such that the difference of any two components of the divisor at infinity is torsion in the Picard group.

\section{Proof Theorem \ref{Main}}

 \paragraph{Three functions on $\Ci$}
Let us consider the functions

\begin{align*}
	f&=\frac{Y+\sqrt{a}X}{W}=\frac{cW}{Y-\sqrt{a}X}\\
	g&=\frac{Z+\sqrt{b}X}{W}=\frac{dW}{Z-\sqrt{b}X}\\
	h&=\frac{\sqrt{b}Y-\sqrt{a}Z}{W}=\frac{(bc-ad)W}{\sqrt{b}Y+\sqrt{a}Z}.
\end{align*}

They are functions on the curve $\Ci$ defined over $\k$.\\
By explicit computation one has

\begin{align*}
	x &=\frac{f^2-c}{2\sqrt{a}f} =\frac{g^2-d}{2\sqrt{b}g} \\ y
&=\frac{f^2+c}{2f} =\frac{bc-ad + h^2}{2\sqrt{b}h} \\ z &=\frac{g^2+d}{2g}
=\frac{bc-ad-h^2}{2\sqrt{a}h},
\end{align*}

so we have that 

\begin{align*}
\k(f)=\k(x,y),\quad \k(g)=\k(x,z),\quad \k(h)=\k(y,z),
\end{align*}
 and $\k(\Ci)$
has degree $2$ over each of them.

 The divisors of the three functions are
\begin{align*}
	(f)&=P_3+P_4-P_1-P_2\\
	(g)&=P_2+P_4-P_1-P_3\\
	(h)&=P_1+P_4-P_2-P_3.
\end{align*}
These three functions are
in $\k[\Ci]^*$, as we can see by their explicit expressions and the expressions for their inverses, or by observing that their divisors are supported on the points at infinity.

By direct computation using the first and then the second definition of $f,g,h$ one immediately checks that they satisfy the linear relations
\begin{subequations}\label{un1}
\begin{align}
	\label{unA}\sqrt{b}f-\sqrt{a}g=h\\
	\label{unB}\frac{c\sqrt{b}}{f}-\frac{d\sqrt{a}}{g}=h.
\end{align}
\end{subequations}

\paragraph{Three functions on $\Ci^{(2)}$}
Let us now consider the Cartesian product $\Ci\times\Ci$ given by equations
\begin{equation*}
	\begin{cases}
		y^2=ax^2+c\\
		z^2=bx^2+d\\
		y'^2=ax'^2+c\\
		z'^2=bx'^2+d.
	\end{cases}
\end{equation*}
and let us indicate by $f',g',h'$ the functions corresponding to $f,g,h$ in
the primed variables.
Let $\Ci^{(2)}$ be the symmetric product of $\Ci$ obtained from $\Ci\times\Ci$
as a quotient by the action of $\Z/2\Z$ that acts swapping the two
coordinates.

The ring of regular functions over $\Ci^{(2)}$ is
$A=\k[\Ci\times\Ci]^{\Z/2\Z}$, the subring of $\k[\Ci\times\Ci]$ consisting of
the functions invariant for this action, and the points of $\Ci^{(2)}$ can be
thought as unordered couples of points on $\Ci$.
\begin{align*}
	\Ci\times\Ci  \overset{\pi}{\longrightarrow}  \Ci^{(2)}\\
	\k[x,y,z,x',y',z']  \supseteq  A.
\end{align*}
Let us denote by $\pi$ the quotient map from $\Ci\times\Ci$ to $\Ci^{(2)}$.

Let $P$ be a quadratic integral point on $\Ci$, and let $P'$ be its
conjugate.

The pair $(P,P')$ is a quadratic integral point on $\Ci\times\Ci$,
and the unordered couple $\{P,P'\}$ is an integral point on $\Ci^{(2)}$ which
is fixed by any Galois automorphism over $\k$ because any such morphism either
fixes both $P,P'$ or swaps them.

The point $\{P,P'\}$ is therefore defined
over $\k$.

\medskip

Let us now define three functions
\begin{align*}
	\alpha&=\frac{cd}{ff'gg'} =&
(y-\sqrt{a}x)(y'-\sqrt{a}x')(z-\sqrt{b}x)(z'-\sqrt{b}x')/c d\\
	\beta&=\frac{c(bc-ad)}{ff'hh'} =&
(y-\sqrt{a}x)(y'-\sqrt{a}x')(\sqrt{b}y+\sqrt{a}z)(\sqrt{b}y'+\sqrt{a}
z')/c(bc-ad)\\
	\gamma&=\frac{d(ad-bc)}{gg'hh'} =&
(z-\sqrt{b}x)(z'-\sqrt{b}x')(\sqrt{b}y+\sqrt{a}z)(\sqrt{b}y'+\sqrt{a}z')/d(ad-bc)\\
    \frac{1}{\alpha}&=\frac{ff'gg'}{cd}=&
(y+\sqrt{a}x)(y'+\sqrt{a}x')(z+\sqrt{b}x)(z'+\sqrt{b}x')/cd\\
    \frac{1}{\beta}&=\frac{ff'hh'}{c(bc-ad)}=
&(y+\sqrt{a}x)(y'+\sqrt{a}x')(\sqrt{b}y-\sqrt{a}z)(\sqrt{b}y'-\sqrt{a}z')/c(bc-ad)\\
    \frac{1}{\gamma}&=\frac{gg'hh'}{d(ad-bc)}=&
(z+\sqrt{b}x)(z'+\sqrt{b}x')(\sqrt{b}y-\sqrt{a}z)(\sqrt{b}y'-\sqrt{a}z')/d(ad-bc).
\end{align*}
We clearly see that they belong to $\k[\Ci^{(2)}]$, and so do their inverses; the functions $\alpha, \beta, \gamma$ are the three regular and nonvanishing functions that we need to follow Vojta's strategy and view $\Ci^{(2)}$ as a subvariety of $\G_m^3$.

The functions $\alpha, \beta$ and $\gamma$ are defined on a surface, so they must be algebraically dependent. Our next step is to find a relation between them.

\medskip

Multiplying together equation \eqref{unA} for the primed and unprimed variables one gets
\begin{equation}\label{troisA}
 hh'=bff'+agg'-\sqrt{ab}(f'g+fg');
\end{equation}
doing the same for equation \eqref{unB} gives
\begin{equation}\label{troisB}
 hh'=\frac{c^2 b}{ff'}+\frac{d^2 a}{gg'}-cd\sqrt{ab}\left(\frac{1}{f'g}+\frac{1}{fg'}\right).
\end{equation}
If we now multiply \eqref{troisB} by $\frac{ff'gg'}{cd}$ and subtract it from \eqref{troisA}, after some tidying up and using the definitions for $\alpha,\beta,\gamma$ we obtain
\begin{equation*}
 \alpha+\beta+\gamma=1.
\end{equation*}

\paragraph{Computing degrees} 
From what we have said until now we have the mappings
\begin{equation*}
	\Ci\times\Ci  \overset{\pi}{\rightarrow}  \Ci^{(2)} 
\xrightarrow{(\alpha,\beta,\gamma)}  H  \subseteq  \G_m^3 
\end{equation*}
where $H$ is the subvariety of $\G_m^3$ defined by $X+Y+Z=1$.

The
corresponding homomorphisms between the rings of regular functions are
\[
	\k\left[X,Y,Z,\frac{1}{XYZ}\right] \xrightarrow[Z\mapsto 1-X-Y]{}
\k\left[X,Y,\frac{1}{XY(1-X-Y)}\right] \xrightarrow[\stackrel{ X\mapsto
\alpha}{ Y\mapsto \beta}]{} A \subseteq \k[\Ci\times\Ci].
\]

To find the degree of the map $\Ci^{(2)}  \xrightarrow{(\alpha,\beta,\gamma)} 
H$ we must find the degree 
\[
	[\k(\Ci^{(2)}):\k(\alpha, \beta,\gamma)].
\]

One finds directly that $(ff')^2=-\frac{c^2\gamma}{\alpha\beta}$, so that $[\k(ff',gg'):\k(\alpha,\beta,\gamma)]=2$, and from
\begin{align*}
 \frac{f^2-c}{\sqrt{a}f}&=\frac{g^2-d}{\sqrt{b}g}\\
\frac{(ff')^2-cf^2}{\sqrt{a}(ff')f}&=\frac{(gg')^2-dg^2}{\sqrt{b}(gg')g}
\end{align*}
follows that $[\k(f,f',g,g'):\k(ff',gg')]\leq 4$, so by our previous remarks we have
\begin{multline*}
   [\k(\Ci\times\Ci):\k(\alpha,\beta,\gamma)]=[\k(f,f',g,g'):\k(\alpha,\beta,
\gamma)]=\\
    =[\k(f,f',g,g'):\k(ff',gg')]\cdot[\k(ff',gg'):\k(\alpha,\beta,\gamma)]\leq 8.
\end{multline*}
Obviously $[\k(\Ci\times\Ci):\k(\Ci^{(2)})]=2$, as there are two ordered pairs for a
generic unordered couple; therefore we can conclude that the degree of the map
given by $\alpha,\beta,\gamma$ between $\Ci^{(2)}$ and $H$ is at most four.

\bigskip

If $P$ is an $S$-integral point on $\Ci^{(2)}$ the values
$\alpha(P),\beta(P),\gamma(P)$ will be $S$-integers, and so will be their
inverses $\frac{1}{\alpha(P)},\frac{1}{\beta(P)},\frac{1}{\gamma(P)}$.

The
point $P$ will then provide a solution in $\O_S^*$ of the equation
\[x_1+x_2+x_3=1.\]
\paragraph{Non degenerate solutions}
Theorem \ref{EvertseSUnitEffective} tells us that there are only finitely-many
non degenerate solutions, and that we can bound their number. If we apply the theorem with $n=3$ we obtain that the number of non degenerate triples $(\alpha,\beta,\gamma)\in H$ is at most $2^{35\cdot 3^4 s}=2^{2835s}$.

We already bounded in the previous paragraph the degree of the map $(\alpha,\beta,\gamma )$, which is at most four. For every integral point on $\Ci^{(2)}$ we have two quadratic integral points on $\Ci$, hence the number of quadratic integral points on $\Ci$ corresponding to non degenerate solutions is at most $2^{2835s+3}$.
\paragraph{Degenerate solutions}
The degenerate solutions to $\alpha+\beta+\gamma=1$ are those with a subsum equal to 0, that is those for which one of the three functions $\alpha, \beta, \gamma$ is equal to $1$; let us then define
\begin{align*}
 W_\alpha: \alpha=1,\qquad W_\beta:\beta=1,\qquad W_\gamma: \gamma=1
\end{align*}
the subsets of $\Ci^{(2)}$ thus obtained.

Using the definition of $\alpha, \beta, \gamma$ we see that, for example,
\begin{align*}
 \alpha&=1, &\beta&=-\gamma \\
 ff'gg'&=cd, &\frac{c(bc-ad)}{ff'hh'}&=-\frac{d(ad-bc)}{gg'hh'}\\
ff'\frac{cgg'}{d}&=c^2, & ff'&=\frac{cgg'}{d}\\
(ff')^2&=c^2, & gg'&=\frac{ff'}{c}d
\end{align*}
and similarly for the other two cases, so we see that $W_\alpha, W_\beta, W_\gamma$ are composed of two subvarieties each, and
\begin{align*}
&W_\alpha =W_x^- \cup W_x^+\\
&W_\beta =W_y^- \cup W_y^+\\
&W_\gamma =W_z^- \cup W_z^+,
\end{align*}
where
\begin{align*}
&W_x^-=\{ff'=c,gg'=d \} &W_x^+&=\{ff'=-c,gg'=-d \}\\
&W_y^-=\{ff'=-c,hh'=ad-bc \} &W_y^+&= \{ff'=c,hh'=bc-ad \}\\
&W_z^-=\{gg'=-d,hh'=bc-ad \}  &W_z^+&=\{gg'=d,hh'=da-bc \}.
\end{align*}
These six subvarieties give all degenerate solutions; to understand them better we use the following simple lemma:
\begin{lemma}
 Let $P=((x,y,z),(x',y',z'))\in \Ci\times \Ci$. Then
\begin{itemize}
 \item $ff'(P)=\pm c \Rightarrow x=\mp x'$ and $y=\pm y'$;
 \item $gg'(P)=\pm d \Rightarrow x=\mp x'$ and $z=\pm z'$;
 \item $hh'(P)=\pm (bc-ad) \Rightarrow z=\mp z'$ and $y=\pm y'$.
\end{itemize}
\end{lemma}
\begin{proof}
 For example, if $ff'=c$ then
\[
 x'=\frac{(f')^2-c}{2\sqrt{a}f'}=\frac{c^2/f^2-c}{2\sqrt{a}c/f}=\frac{c-f^2}{2\sqrt{a}f}=-x,
\]
and similarly for the other five cases.
\end{proof}
If we define six subvarieties of $\Ci\times \Ci$
\begin{align*}
 V_x^\pm &: x=\pm x',y=\mp y', z=\mp z'\\
 V_y^\pm &: x=\mp x',y=\pm y', z=\mp z'\\
 V_z^\pm &: x=\mp x',y=\mp y', z=\pm z',
\end{align*}
then the lemma tells us that each of them corresponds through $\pi$ to the similarly named subvariety $W_{x,y,z}^\pm$ of $\Ci^{(2)}$ (we denote by $V_{x,y,z}^\pm$ any of the six varieties $V_x^-,V_x^+,V_y^-,V_y^+,V_z^-,V_z^+$, and we do the same with $W_{x,y,z}^\pm$).
In short we have that
\begin{equation*}
 \pi|_{V_{x,y,z}^\pm}:V_{x,y,z}^\pm \to W_{x,y,z}^\pm
\end{equation*}
Note also that $V_{x,y,z}^\pm\simeq \Ci$ through the projection on the first component of $\Ci\times \Ci$, so that these six curves all have genus one.

\paragraph{Points on $W_{x,y,z}^+$}
Consider for example
\begin{align*}
  \pi|_{V_{x}^+}:V_{x}^+ &\to W_{x}^+\\
  ((x,y,z),(x,-y,-z)) &\mapsto \{(x,y,z),(x,-y,-z) \}.
\end{align*}
A generic point in $W_{x}^+$ has two preimages, obtained by exchanging the order of the pair; the points which have just one primage are those such that $y=z=0$, but we see from the defining equations \eqref{DefiningEquation} of $\Ci$ that this can never happen because $by^2-az^2=bc-ad\neq 0$.
The map $\pi|_{V_{x}^+}$ is therefore an unramified covering of $W_{x}^+$, and the Riemann-Hurwitz formula tells us that
\begin{align*}
 0=2g(V_{x}^+)-2=2(2g(W_{x}^+)-2)\\
g(W_{x}^+)=1.
\end{align*}
The same is true for $W_{y}^+$ and $W_{z}^+$, because there is no point on $\Ci$ where two of $x,y,z$ both vanish.

We have thus shown that the three subvarieties $W_{x,y,z}^+$ are all curves of genus 1, so they carry only a finite number of integral points, which in turn correspond to a finite number of quadratic integral points on $\Ci$.

\medskip

We might also argue the finiteness of quadratic integral points coming from the curves $W_{x,y,z}^+$ as follows: 
let $K=\k(\sqrt{\epsilon})$ be a quadratic extension of $\k$, we may suppose
$\epsilon\in \O_S$; let $P=(x_0,y_0,z_0)$ be a $K$-integral point on $\Ci$, and
$P'=(x_0',y_0',z_0')$ its conjugate.
Suppose that $(P,P')$ belongs to $V_x^+$; if it is so, after enlarging $S$ to an $S'$ so that $\O_{S'}$ has trivial class group, we can write
$x_0=t, y_0=u\sqrt{\epsilon}, z_0=v\sqrt{\epsilon}$ for some $t,u,v\in\O_{S'}$.

Substituting back into the equations for $\Ci$ we get
$b\epsilon u^2 -a \epsilon v^2=bc-ad$, so $\epsilon$ divides $bc-ad$; different $\epsilon$ differing only by a square factor give the same extension, and given that $\O_{S'}^*/(\O_{S'}^*)^2$ is finite we see that $K$ must belong to a finite set of
extensions of $\k$, whose cardinality may be bounded in terms of the cardinality of $S'$.

We know, by Siegel's theorem, that there are only finitely many integral
points on $\Ci$ defined over a fixed number field; since $\Ci$ has four points at infinity, this number may be bounded effectively, as done in \cite{CZ2}, in terms of the degree of $K$ and the cardinality of the extension of $S'$ to $K$, which in turn are both bounded in terms of $\k$ and $|S'|$.

The computations involved in this bound are quite heavy, but again the special structure of the curve $\Ci$ provides us with a simpler argument: equation \eqref{unA} gives the unit equation
\begin{equation}\label{aaa}
\frac{\sqrt{b}f}{h}(P)-\frac{\sqrt{a}g}{h}(P)=1,
\end{equation}
so the number of solutions may be easily bounded using again Theorem \ref{EvertseSUnitEffective}.
For a fixed $K$ and extension $S''$ of the set of absolute values $S'$, there are at most $2^{35\cdot 2^4 |S''|}=2^{1120|S'|}$ solutions; the number of extensions $K$ is at most $|\O_{S'}^*/(\O_{S'}^*)^2|=2^{|S'|}$; the functions $f/h$ and $g/h$ have degree two, as we can see computing their divisors, so each solution gives at most two points. Combining all, we have that the number of quadratic integral points $P$ such that $(P,P')$ lies
on any of the $W_{x,y,z}^+$ is at most $3\cdot 2^{1121|S'|+1}=3\cdot 2^{1121(s+H-1)+1}$, where $H$ is the class number of $\O_S$.

\medskip

I thank the anonymous referee for suggestions about this bound.

\medskip

We should also remark that, since the relevant curve has genus 1, it is in fact possible, as is well known, to bound {\em the height} of the solutions themselves; tha same can be done on equation \eqref{aaa}, using effective results on {\em two-term} unit equations (derived from Baker's theory of linear forms in logarithms). So the quadratic integral points arising from $W_{x,y,z}^+$ are, in fact, effectively computable (in contrast with those arising from the ``non degenerate solutions'').

\paragraph{Points on $W_{x,y,z}^-$}
Reasoning as we did before, we consider, for example
\begin{align*}
  \pi|_{V_{x}^-}:V_{x}^- &\to W_{x}^-\\
  ((x,y,z),(-x,y,z)) &\mapsto \{(x,y,z),(-x,y,z) \}.
\end{align*}
This is again a map of degree $2$.

In this case however, we see that it is ramified at points where $x=0$, and there are $4$ such points on $\Ci$, namely $(0,\pm\sqrt{c},\pm\sqrt{d})$, each of them of course ramified of index 2. 
Therefore this time applying the Riemann-Hurwitz formula we obtain
\begin{align*}
 0=2g(V_{x}^-)-2=2(2g(W_{x}^--2)+4\\
g(W_{x}^-)=0;
\end{align*}
as $W_{x}^-$ has genus 0, it may
contain infinitely many integral points.

The composition of the maps
\begin{align*}
 \Ci \to V_x^-  \to W_x^- \xrightarrow{\pi|_{W_x^-}} \{by^2-az^2=bc-ad\}\simeq\P_1\\
 (x,y,z)  \mapsto ((x,y,z),(-x,y,z))  \mapsto \{ (x,y,z),(-x,y,z)\}\mapsto (y,z)
\end{align*}
gives a map of degree two from $\Ci$ to $\P_1$ that takes quadratic integral points arising from $W_x^-$ to integral points in $\P_1$; this map together with the same compositions for $W_y^-$ and $W_z^-$, that is,
\begin{align*}
 \Ci\to \{z^2=bx^2+d\}\simeq\P_1\\
(x,y,z)\mapsto (x,z),
\end{align*}
and
\begin{align*}
 \Ci\to \{y^2=ax^2+c\}\simeq\P_1\\
(x,y,z)\mapsto (x,y),
\end{align*}
are the maps in Theorem \ref{ThmCZQuadratic}.

\paragraph{Acknowledgments} The author thanks Professors Corvaja, Evertse and Zannier for the substantial advice received in the preparation of this article, and the anonymous referee for the thorough report and for having suggested some improvements.

\bibliographystyle{amsalpha}
\bibliography{Double_Pell-art}
\end{document}